\documentclass[12pt]{article}
\usepackage{stackrel}
\usepackage{amsmath}
\usepackage{amsfonts}
\usepackage{amssymb}
\usepackage{amsthm}
\usepackage{color}

\newtheorem{exm}{Example}[section]

\newcommand{\R}{{\rm I}\kern-0.18em{\rm R}}
\newcommand{\1}{{\rm 1}\kern-0.25em{\rm I}}
\newcommand{\E}{{\rm I}\kern-0.18em{\rm E}}
\newcommand{\p}{{\rm I}\kern-0.18em{\rm P}}


\usepackage{epsfig}

\usepackage{fancyhdr}

\usepackage{graphics}

\usepackage{amsopn}  %%% for new Math operators
\usepackage{amssymb}
\usepackage{amsxtra}
\usepackage{amsfonts}

\usepackage{afterpage}
\usepackage{eucal}
\usepackage{eufrak}

\usepackage{enumerate}

\allowbreak

\def\fnote#1{\footnote}

\newcommand{\bea}{\begin{eqnarray}}

\newcommand{\eea}{\end{eqnarray}}

\newcommand{\beas}{\begin{eqnarray*}}
\newcommand{\eeas}{\end{eqnarray*}}

\title{Tempered distributions: does universal tempering procedure exist?}
\author{Lev B. Klebanov\footnote{Department of Probability and Statistics,
Charles University}, Lenka Sl\'{a}mov\'{a}\footnote{Department of Probability and Statistics,
Charles University}
}
\date{}
\begin{document}
\maketitle

\begin{abstract}
Since the turn of the century, there has been increased interest in the application of heavy-tailed distributions, particularly stable distributions, to problems in physics and finance.  Although, the tails of stable distributions provide a better fit to real-world data, they are  too fat to describe empirical distributions. To remedy this drawback of stable distributions, so-called tempered variants of stable distributions have been proposed. In this paper, we argue that the tempering should be connected to the model leading to heavy tailed distribution and propose several tempering procedures in connection with corresponding models
\end{abstract}


\section{Introduction}
\setcounter{equation}{0}

Stable distributed are widely used in many branches of the sciences and economics because of the overwhelming empirical evidence that real-world distributions exhibit fatter tails than that of the Gaussian distribution.  Examples in financial economics include the works of Mandelbrot \cite{M1963}, Fama \cite{F1965}, Samuelson \cite{S1967},  Embrechts et al. \cite{Em1997}, Rachev and Mittnik \cite{RM2000}. However, it is clear, that all amount of money in the world is finite. This implies paradoxically that heavy tailed distributions may provide a good model for large, but not too large money sums. Accordingly, one may expect, that the tails of stable distributions are too heavy for consideration of large money amounts. Connection of this idea with applicability of generalized limit theorem was proposed in \cite{KRSz, KRSa}. Despite the heavier tails than the Gaussian distribution, stable distributions have been found too heavy than found in real-world distributions (see, for example, \cite{CGMY}). The most popular approach for correcting for heavier tails than warranted is the so-called tempering procedure that involves modifying the characteristic function of stable distribution by changing of L\'{e}vy measure in  L\'{e}vy-Khinchine representation of stable distributions. The characteristic function $\phi_{CTS}$ for a tempered stable distribution is given by 
\[
\phi_{CTS}(u) = \exp\Bigl( iu\mu +C_1\Gamma(-\alpha)((\lambda_+-iu)^\alpha-\lambda_+^\alpha) 
  +C_2\Gamma(-\alpha)((\lambda_-+iu)^\alpha-\lambda_-^\alpha) \Bigr),
 \]
for some $\mu\in\R$. Moreover, $\phi_{CTS}$ can be extended to the region 
\[ \{z\in\mathbb{C}: \;\operatorname{Im}(z)\in(-\lambda_-,\lambda_+)\}. \]
Despite the popularity of this approach, we demonstrate that this procedure is  not universal and any modification of the stable distribution should be made on an individual basis. We also are that for non-stable distributions with heavy tails, the tempering procedure makes no sense for such distributions.

\section{Examples of models leading to heavy tailed distributions}
\setcounter{equation}{0}

Here we give examples of toy-models (i.e. small models ignoring many details of the process studied, but taking into account one of the most essential part of it) for probability distributions with heavy tails. 

\begin{exm}\label{ex3} \textbf{The First Passage Time Distribution}. Here we are interested in the probability that a particle reaches a given point $c$ at a time $t$. It is well-known that first-passage time for a Brownian particle follows a L\'{e}vy distribution (see \cite{Fel2}). This distribution is $\frac{1}{2}$-stable. It is concentrated on the positive semi-axis and has characteristic function
\begin{equation}\label{eq3}
f(t) =\exp\{ -\sqrt{-2bit}\}, \;\; t \in \R^1, \;\; b>0.
\end{equation}
This distribution has a heavy tail of order $1/2$.
\end{exm}

\begin{exm}\label{ex4} \textbf{Discrete variant of Example \ref{ex3}}. Now consider the random walk on real line with equal probabilities of moving to the right or to the left. We are interested in the distribution of the time of the first passage through 1. It is known that the probability generating function of this distribution is 
\begin{equation}\label{eq4}
{\mathcal P}(z)=\frac{1-\sqrt{1-z^2}}{z}.
\end{equation}
Clearly, the corresponding distribution has a heavy tail of order $1/2$.
\end{exm}

\begin{exm}\label{ex5} \textbf{Sub-Gaussian distributions}. Let $X$ be a Gaussian random variable with zero mean. Suppose that $A$ is $\alpha$-stable positive random variable (it is necessary to have $\alpha \in (0,1)$, and skewness parameter $\beta=1$). If $X$ and $A$ are independent, then the product $A^{1/2}X$ has a symmetric $2\alpha$-stable distribution. This is the so-called sub-Gaussian distribution. Clearly, the sub-Gaussian distribution can be represented as a product of two stable distributed random variables. Such representation is not unique.
\end{exm}

\begin{exm}\label{ex2}\textbf{LePage series}.
The LePage series allows to obtain many examples of toy-models leading to heavy tailed distributions
(see, for example, \cite{ST}). Suppose, we have a Poisson process, and its arrival times are $\Gamma_j$, $j=1,2, \ldots$. Let $\{X_1, X_2, \ldots \}$ be a sequence of independent identically distributed (i.i.d.) random variables and suppose that this sequence is independent with Poisson process mentioned above. Consider the following series
\begin{equation}\label{eq2}
\sum_{j=1}^{\infty}\frac{X_j}{\Gamma_j^{1/\alpha}},
\end{equation}
where $\alpha$ is a real number between 0 and 2 (strictly), and random variable $X_1$ has an absolute moment of order greater that $\alpha$. If the series (\ref{eq2}) converges, then its sum follows a strictly stable distribution with index of stability equals to $\alpha$.
\end{exm}

The sum of the series (\ref{eq2}) allows for the following three interesting interpretations. 
\begin{enumerate}
\item Consider electrical charges $X_j$, $(j=1,2, \ldots)$ on a straight line, arranged with accordance to the distances $\Gamma_j$  from the origin, where is a fixed charge $Q$. According to Coulomb's law, the force acting on the charge $Q$ is
\[ Y=\sum_{j=1}^{\infty}\frac{X_j}{\Gamma_j^{2}}.\]
Under condition that all charges $X_j$ are chosen as i.i.d. random variables, we obtain that the random force $Y$ has a strictly stable distribution with the parameter $\alpha =1/2$.
\item Similar to the first interpretation suppose that instead of charges  we have masses, and use Newton law instead that of Coulomb's law. Again, the force acting on the mass $Q$ will have a strictly stable distribution with stability parameter $\alpha = 1/2$. In this situation, all $X_j$ are positive, and, therefore, random variable $Y$ is positive, too. Now we may conclude, that $Y$ has L\'{e}vy distribution. Note that we may consider the Poisson field on Euclidean plane or in Euclidean space. The distances $\Gamma_j$ from a particle to the origin can be considered as the arrival times of a Poisson process. Therefore, we will have the same interpretations for the mass (or charges) distributed according to the Poisson field.
\item Consider Euclidean plain $\R^2$, and suppose that there is a mobile phone base station at the origin. Suppose further that mobile phones are distributed on the plane according to a Poisson field. The $j$th phone sends  signal $X_j$ to the base station. The signal decays inversely to some exponent of the distance between the phone and base station (i.e., the signal, coming to the station from a phone is $X_j/\Gamma_{j}^{1/\alpha})$. From technical experiments it is known, that $1/\alpha$ is a number around 2.6. Now we can say that the summand signal coming to the base station has strictly stable distribution.
\end{enumerate}

\begin{exm}\label{ex1} \textbf{Pareto distribution}. Let $S=\{X_1, X_2, \ldots \}$ be a sequence of i.i.d. positive random variables. Suppose also that $\{ \nu_p, \; p\in (0,1) \}$ is a family of random variables, independent of the sequence $S$ and having geometric distribution with parameter $p$: $\p\{\nu_p=k\} = p(1-p)^{k-1}$, $k=1,2, \ldots$. Consider the normalized product of random number of the elements of $S$, that is,
\begin{equation}\label{eq1}
Z_p=\prod_{j=1}^{\nu_p}X_j^{p}. 
\end{equation}
It is known (see \cite{KMR}) that the limit distribution of $Z_p$ as $p \to \infty$ is a Pareto distribution with shape parameter $a=1/\gamma$, where $\gamma=\E \log X_1$, assuming this parameter is positive, that is,
\begin{equation}\label{eq1a}
\lim_{p \to 0}\p\{Z_p<x\}=\begin{cases} 1-x^{-a} & \text{for}\; x>1,\\ 0 & \text{for}\; x<1, \end{cases} 
\end{equation}
in the case when $\gamma >0$.
This particular case is interesting for us because it leads to heavy tailed distribution. 

Also note that if random variable $X_1$ has a Pareto distribution with parameter $a>0$, then
\[ X_1 \stackrel{d}{=} \prod_{j=1}^{\nu_p}X_j^{p}, \]
where $X_1, X_2, \ldots $ are i.i.d. random variables, and $\{\nu_p, \; p\in (0,1)\}$ is a family described above. The sign $\stackrel{d}{=}$ is used to denote the equality in distribution.
\end{exm}

Note that the distribution of the value of $Z_p$ represents the distribution of capital, obtained in random market over $\nu_p$ units of time. The mean value of $\nu_p$ is large in this situation (or, equivalently, the parameter $p$ is small). On economical interpretation of $Z_p$ and its limit distribution the reader is addressed to \cite{KMR}

Example \ref{ex1} shows that the reason for the distribution having heavy tails is the randomness of the number of multipliers in the product (\ref{eq1}). This number can be arbitrarily large with a positive probability, while the multipliers may have thin tails. However, note that we do not have  a sum of a random number of random variables, but a product of them. In the case of sums it is impossible to obtain heavy tailed limit distributions using summands with thin tails. Accordingly  {\it the limit theorem on convergence to stable distribution cannot be considered as an explanation for the appearance  of heavy tailed distributions in applied problems}. The central limit theorem shows, we need to have the summands with heavy tails for sum to be convergent to stable distribution.

\begin{exm}\label{exB} \textbf{Selling-short-strategy in financial market.} Suppose that an investor follows a selling-short-strategy. This means that the investor has borrowed some assets and sells it at Price of $P_{*}$  with the expectation that the price will decline below. Should the price decline, the investor will buy the assets at the lower price and returns the asset to the original owner. Let us study this strategy in more detail for the discrete case. Suppose that $\{P_j, \; j=1, 2, \ldots\}$ is a sequence of prices for an asset at times $1,2, \ldots$. In our toy-model we consider $P_1, P_2, \ldots $ as a sequence of i.i.d. random variables (which is, of course, a simplification). The short-seller borrows the asset, and will sell it at time $1,2, \ldots $ prior to the first appearance of the event $\{P<P_{*}\}$. The moment of this appearance will be denoted by $\nu+1$.  It is clear that the time elapsed before the onset of the event may be considered a random variable. A first approximation for describing this random variable is to consider it as a number of experiments in a Bernoulli scheme until the first appearance of a success. The probability of success is $p=\p\{P<P_{*}\}$. In other words, the random number of experiments has geometric distribution (starting at 1) with parameter $p \in (0,1)$. Now we see that the sum obtained by the short-seller in this strategy is 
\begin{equation}\label{eqBm}
S=\sum_{j=1}^{\nu_p} P_j X_j.
\end{equation}
Here $P_j$ is the price at the time $j$, $X_j$ is the quantity of assets sold by the short-seller at this time, and $\nu_p$ is a random variable with geometric distribution $\p\{\nu_p = k\} = p(1-p)^{k-1},$ $k=1,2, \ldots$. For simplicity, we suppose that $\nu_p$ is independent on both sequences $\{P_j, \; j=1,2, \ldots\}$ and $\{X_j, \; j=1,2, \ldots\}$, and last two sequences are independent on each other as well. Let us make some suppositions on the character of the sequence  $\{X_j, \; j=1,2, \ldots\}$. Random variable $X_j$ is an integer equals to the number of assets sold at moment $j$. It is more or less natural to consider $X_j$ as a number of successes in Bernoulli scheme, but with the probability of fail depending on the number of the experiment. Namely, we consider the process of selling at the time $j$ in the following way. It is a sell of the first asset, after that a sell of the second asset, and so on, until nobody will buy assets. However, the probability to sell second asset is supposed to be higher, than for the first asset; probability to sell third asset is higher than for the second one; and so on. It can be connected to the fact that most gamblers (according to some sources, up to 80\%) prefer strategy ``to buy". For the first approximation, we suppose, that the probability to fail while selling $k$th asset is $\gamma/k$, where $\gamma$ is the probability to fail while selling the first asset. Under these suppositions, we see that all random variables $X_{j}, \; j=1,2, \ldots$ are i.i.d. with the distribution of the form
\[  \p\{ X_j =k \} = \frac{\gamma}{k} \prod_{i=1}^{k-1}(1-\frac{\gamma}{i}). \]
It is clear that probability generating function of $X_j$ is
\begin{equation}\label{eqSib}
Q(z) = 1-(1-z)^{\gamma}.
\end{equation}
Let $L_{P}(s)$ be Laplace transformation of random variable $P_1$. Then Laplace transform of the product $P_jX_j$ is
\[ L_{PX}(s)= \E \exp\{ -s P_j X_j\} = \sum_{k=1}^{\infty}L_{P}(sk)\p\{X_j =k\} =\sum_{k=1}^{\infty}L_{P}(sk) \frac{\gamma}{k} \prod_{i=1}^{k-1}(1-\frac{\gamma}{i}). \]
Now it is easy to calculate Laplace transform of the sum (\ref{eqBm})
\begin{equation}\label{eqS}
L_{S}(s)=\E \exp\{ -sS\}=\frac{p L_{PX}(s)}{1-(1-p)L_{PX}(s)}.
\end{equation}
Let us find transforms of $L_{PX}(s)$ and $L_S(s)$ for a particular case of distribution of prices $P_j$, namely, for the case of exponentially distributed prices. In this case we have
\[  L_{P}(s)=\frac{1}{1+as},\]
where $a>0$ is a parameter of exponential distribution. The expression for Laplace transform of $P_jX_j$ is
\[  L_{PX}=\sum_{k=1}^{\infty}\frac{1}{1+ask} \frac{\gamma}{k} \prod_{i=1}^{k-1}(1-\frac{\gamma}{i}) =1-\frac{\Gamma (1+1/(as))\Gamma (1+\gamma)}{\Gamma(1+\gamma+1/(as))},\]
where $\Gamma (z)$ is Euler gamma-function. From (\ref{eqS}) we get
\begin{equation}\label{eqBF}
L_S(s) = \frac{p \Bigl( \Gamma(1+\gamma+1/(as))-\Gamma(1+\gamma)\Gamma(1+1/(as))\Bigr)}{p\Gamma(1+\gamma+1/(as))+(1-p)\Gamma(1+\gamma)\Gamma(1+1/(as))},
\end{equation}
where the parameters satisfy to the conditions: $a>0$, $0<\gamma<1$. It is not difficult to calculate the limit
\[ \lim_{s \to 0}\bigl( 1-L_S(s)\bigr)/s^{\gamma} = a^{\gamma}\Gamma (1+\gamma)/p. \]
In view of this, the distribution of random variable $S$ has heavy (power) tail of order $\gamma$.

Random variable $S$ is not the profit of the short-seller. This is the amount he/she received from the sale of assets. Short-seller has to buy the assets to give them back. It is necessary to buy $\sum_{j=1}^{\nu}X_j$ assets for the price less or equal to $P_{*}$. The profit will be greater or equal to $\sum_{j=1}^{\nu}(P_j-P_{*}) X_j$. Under conditions imposed above, this is a random variable with heavy tailed distribution.

\end{exm}

The list of toy-models may be prolonged, however, it is enough to explain multiple ways for tempering such heavy tails distributions. Some of the ways will be proposed in the next section.

\section{Tempering procedures}
\setcounter{equation}{0}

In this section we propose different variants of tempering procedures. The procedures are connected with an attempt to make underlying toy-model more precise in a sense.

\subsection{Considerations of Example \ref{ex3}}
Example \ref{ex3} has physical sense, and clearly,  there are some idealizations in the Example. One of idealizations consists in supposition that Brownian motion has no drift. The presence of heavy tails in the first passage time distribution is connected exactly with this fact. While considering the motion with drift, we obtain Inverse Gaussian Distribution as the first passage time distribution (see, for example \cite{Fel2}). The form of probability density of this distribution shows the presence of exponential tail instead of heavy tail of L\'{e}vy distribution. The density of Inverse Gaussian Distribution has form
\begin{equation}\label{eqIG}
p_{IG}(x,\lambda ,\mu)=\sqrt{\frac{\lambda}{2\pi x^3}}\exp\Bigl(-\frac{\lambda (x-\mu)^2}{2x\mu^2}\Bigr)
\end{equation}
for $x>0$. The density of  L\'{e}vy distribution is
\begin{equation}\label{eqLevy}
p_{L}(x,\sigma)=\sqrt{\frac{\sigma}{2\pi x^3}}\exp \Bigl(-\frac{\sigma}{2x}\Bigr),
\end{equation}
$x>0$. For the case $\lambda = \sigma$ we have
\[ p_{IG}(x,\sigma, \mu)=p_{L}(x,\sigma)\exp\{\sigma/\mu\}exp\{-\sigma x /(2\mu)\}, \]
where the term $exp\{-\sigma x /(2\mu)\}$ responds for the drift, and $\exp\{\sigma/\mu\}$ is a normalizing constant. In terms of characteristic functions, this transformation corresponds to transition from real axis to a line, parallel to it  in complex plane. This can be seen from their characteristic functions:
\[ f_{IG}(t) = \exp\Bigl( \frac{\sigma (1-\sqrt{1-2it\mu^2/\sigma}}{\mu} \Bigr) \]
for Inverse Gaussian Distribution ($\lambda = \sigma$), and
\[ f_{L}(t) =\exp \bigl(-\sqrt{-2 i \sigma t}  \bigr)  \]
for  L\'{e}vy distribution. Transformation from $f_{L}$ to $f_{IG}(t)$ is exactly ``classical" tempering procedure (see \cite{BL, CGMY, KRBF} ).  Note, that in the considered situation the transformation looks absolutely natural, and we have no objection about it.

In a similar way may be considered the case of positive stable distributions with index of stability $\alpha \in (0,1)$. Namely, let us consider stable random variables with Laplace transform
\[ L_{st}(s) = \exp(-A s^\alpha) .\]
Multiplying the probability density of this distribution by $\exp (-a x)$ and making renormalization we obtain new probability density with Laplace transform
\[ L_{tst}(s) = \exp(-A (s+a)^{\alpha}) \exp(A a^{\alpha}).\]
Corresponding distribution has exponential tail. It is natural to call this distribution ``tempered" stable distribution. However, this ``tempering" has no such nice physical interpretation as in the case of $\alpha =1/2$, and the question on uniqueness and natural character of this procedure is open. 

\subsection{Considerations of Example \ref{ex4}}

As it has been mentioned above, Example \ref{ex4} is discrete variant of Example \ref{ex3}. Obviously, the distribution with probability generating function
\begin{equation}\label{eqDL}
{\mathcal P}(z)=\frac{1-\sqrt{1-z^2}}{z}
 \end{equation} 
has heavy tail. Because this distribution has the first passage time interpretation, it is natural to propose tempering procedure by introducing a drift into random walk. To this aim, it is sufficient to change the distribution by the first passage time for the case of non-equal probabilities of moving to the right and to the left. For this case, it is known (see \cite{Fel2}) that the probability generating function has the form
\begin{equation}\label{eqDIG}
{\mathcal Q}(z)=\frac{1-\sqrt{1-4p(1-p)z^2}}{2(1-p)z}.
\end{equation}
Here $p$ is the probability of moving to the right. We suppose, $1/2 < p <1$. In terms of characteristic functions we have
\begin{equation}\label{eqDLcf}
f(t)=\frac{1-\sqrt{1-e^{2it}}}{e^{it}}
\end{equation}
for characteristic function corresponding to (\ref{eqDL}), and
\begin{equation}\label{eqDIGcf}
g(t)= \sqrt{\frac{p}{1-p}}\,\cdot \frac{1-\sqrt{1-\exp(2i(t-ia))}}{\exp(i(t-ia)},
\end{equation}
where $a=(1/2)\log(4p(1-p))$, for the case of (\ref{eqDIG}).

We see, that the procedure used in Example \ref{ex3} is applied here. In other words, classical tempering works nicely in this Example \ref{ex4} as well.

However, the situation with tempering in Example \ref{ex4} depends on the model interpretation. Namely, it may be imagined, that the player (in game interpretation of the random walk) has only bounded sum of money to pay for each move. The drift is absent. Suppose, the sum of money is $M>0$, and the player has to pay one unit for each move. Suppose also there is no drift in the game (walk). So, we have truncated first passage time distribution. Probability generating function now has the form
\begin{equation}\label{eqTrDL}
{\mathcal P}_M(z)=\sum_{k=1}^{[M/2-1]} (-1)^{k+1} \binom{\frac{1}{2}}{k}z^{2k-1}+\Bigl( 1-\sum_{k=1}^{[M/2]-1} (-1)^{k+1} \binom{\frac{1}{2}}{k} \Bigr) z^{2[M/2]-1}.
\end{equation}
It is obvious, this distribution is not classical tempered variant of that from Example \ref{ex4}. 

Our opinion is that the variant of tempering should be chosen in accordance to the needs of the model correction. 
The second variant of tempering is nor better neither worth than the first one, if there is no information about preferable type of correction. Now one can understand, such problem may arise also while constructing tempering procedure in Example \ref{ex3}. For example, we may have time limitation for a particle passes through given level, what  leads to truncation.

\subsection{Considerations of Example \ref{ex5}}

In Example \ref{ex5} we have a product $XA^{1/2}$, where $X$ is zero-mean Gaussian random variable, and $A$ is positive stable distribution. Therefore, it is natural to suppose, that tempering procedure for the product has to be based on a transformation of stable random variable $A$. As it was mentioned above, we have two variant for the tempering procedure. First one is a classical tempering procedure, while the second is a truncation of the distribution of $A$. Third possibility is just a truncation of the distribution of the product $XA^{1/2}$. Below we show that all three possibilities leads to different results.

Let us consider these possibilities in details. Suppose that a positive random variable $A$ has stable distribution with the index of stability $\alpha \in (0,1)$ and Laplace transform $L(s)=\exp\{-s^{\alpha}\}$. Suppose that a random variable $X$ has standard Gaussian distribution. Then the product $X\cdot A^{1/2}$ has symmetric stable distribution with index of stability $2 \alpha$ and characteristic function
\[ f(t) = \exp \{ -t^{2\alpha}/2^{\alpha}\}. \]
\begin{enumerate}
\item \label{it1} Let start with first tempering procedure of the random variable $A$. Laplace transform of tempered random variable $A_a$ has form
\[ L_{a}(s) = \exp\{-(s+a)^{\alpha}\} \exp\{ a^{\alpha}\},\]
where $a$ is a parameter of tempering procedure. Characteristic function of the product $X A_{a}^{1/2}$ is 
\begin{equation}\label{eq1tem}
f_a(t) = \exp \{-(t^2/2+a)^{\alpha}\}\exp\{a^{\alpha}\}.
\end{equation}
The distribution obtained has exponential tails. It converges to initial symmetric $\alpha$-stable distribution as $a \to 0$. This variant of tempering is of an interest because it proposes an application of the procedure to one (stable) multiplier, not to the whole product.
\item\label{it2} The second possibility is to rewrite the product $X A^{1/2}$ in the form $Y B^{1/\beta}$, where $Y$ is a symmetric stable random variable with the index $\beta \in (0,2)$ and $B$ is a positive stable random variable with index $\gamma \in (0,1)$. In the case when $2 \alpha = \beta \gamma$ it is possible to choose the parameters of $Y$ and $B$ so that $X A^{1/2}\stackrel{d}{=}Y B^{1/\beta}$. Proposed tempering procedure consists now in tempering of $B$ in the way, proposed in item \ref{it1}. Obtained product has the tail equivalent to that of the random variable $Y$. Tempering here does not provide exponential tails. Tails remains to be heavy, but not so heavy like for initial distribution of $X A^{1/2}$.
\item\label{it3} Third procedure consists in truncation of the random variable $A$, i.e. instead of $A$ we consider $A_M=\min (A,M)$ for a a fixed constant $M>0$. Characteristic function of the distribution of $X A_M^{1/2}$ is 
\[  \E \exp\{ itX A_M^{1/2}\} = \E \exp\{ -t^2/2 A_M^{1/2}\}= \int_{0}^{M}\exp\{- x t^2/2\} dF_{\alpha}(x),   \]
where $F_{\alpha}(x)$ is distribution function of $A$. It is clear, that the distribution of $A_M$ has Gaussian tails. Corresponding distribution is {\it not} infinitely divisible.
\end{enumerate}
Some other variants of tempering are possible as well, but we do not discuss them here.

\subsection{Considerations of Example \ref{ex1}}
It is more or less clear, that the limit distribution of $Z_p$ has heavy tails because random number of multipliers $\nu_p$ is large in mean, and $\nu_p \to \infty$ in probability as $n \to \infty$. Of course, this situation represents an idealization of the reality. It is natural to suppose, that random variable $\nu_p$ is bounded with probability 1 by a constant. Therefore, it seems to be natural to change geometric random variable $\nu_p$ by bounded random variable $\mu_{p,M}$:
\begin{equation}\label{eqMu}
\p\{ \mu_{p,M}=k\} = \frac{p(1-p)^{k-1}}{1-(1-p)^M}, \;\;( k=1, \ldots , M) 
\end{equation}    
for $p \in (0,1)$ and integer $M>1$. Corresponding probability generating function has the form
\[ {\mathcal P}_{p,M}(z) =\frac{pz(1-(1-p)^M z^M)}{(1-(1-p)^M)(1-(1-p)z)}. \]
It is clear, that 
\begin{equation}\label{eq5}
{\mathcal P}_{p,M}(z) \stackrel[M \to \infty]{}{\longrightarrow} \frac{pz}{1-(1-p)z}, 
\end{equation}
i.e. the limit distribution in this case is geometric one. However,
\begin{equation}\label{eq6}
{\mathcal P}_{p,M}(z) \stackrel[p \to 0]{}{\longrightarrow} \frac{\sum_{k=1}^M z^k}{M},
\end{equation}
and the limit distribution is uniform over integers from 1 to $M$. 

From (\ref{eq5}) and (\ref{eq6}) is clear, that it is impossible to change the order of limits over $M$ and over $p$. What is a tempered variant of Pareto distribution? Unfortunately, there is no definite answer to this question. It is possible to choose many different possibilities. For example, one can take a sum of random number (\ref{eqMu}) of arbitrary summands. Because the summands may have different distributions, there are infinite many possibilities for the tempering definitions. Another way is the following. One may consider Mellin transform of the Pareto distribution as a characteristic function of an exponential distribution. This distribution is geometric stable, and one can make classical geometric tempering procedure to get corresponding analogue of exponential distribution. Now it is necessarily to come back by considering the characteristic function as Mellin transform of tempered Pareto distribution. Of course, there are many other possibilities to make tempering procedure for Example \ref{ex1}.

\subsection{Considerations of Example \ref{exB}}

It is clear that the reason for the distribution with Laplace transform (\ref{eqBF}) to have heavy tails is that the distribution of random variables $X_j$, $(j=1,2, \ldots)$ (it is so-called Sibuya distribution) has such tails. The most natural way to introduce tempering procedure for this example is to change Sibuya distribution by its truncated variant, i.e. instead of probability generating function (\ref{eqSib}) use 
\begin{equation}\label{eqTrB} 
Q_M(z)= \frac{\sum_{k=1}^{M} (-1)^{k+1}\binom{\gamma}{k}z^k}{\sum_{k=1}^{M} (-1)^{k+1}\binom{\gamma}{k}}.
\end{equation}
Of course, corresponding distribution of $S$ has exponential tails.

As usual, this is not unique way of tempering. We may choose another probability generating function
\begin{equation}\label{eqTeB}
R_a(z) = \frac{1-(1-az)^{\gamma}}{1-(1-a)^{\gamma}},
\end{equation}
where $a>0$ is tempering parameter. 

Variant (\ref{eqTeB}) is more similar to classical tempering procedure. However, the variant (\ref{eqTrB}) seems to be more natural one, because short-seller may have only a bounded number of the assets. The number of buyers is bounden as well. Therefore, they can buy bounded number of assets.

\subsection{Other ways of tempering}

Classical procedure of tempering is based on exponential truncation of L\'{e}vy measure in L\'{e}vy-Khinchine representation of stable characteristic function. Obviously, it is possible to make some other modifications of this measure. For example, it may be changed by a discrete one. Such changes were proposed in \cite{SK} in connections with definitions of discrete stable distributions. Unfortunately, in \cite{SK} there are no toy-models explaining such way of tempering.

\section{Conclusion} 
\setcounter{equation}{0}
Morality of the article is that there exist no universal tempering procedure. There are many such procedures for a distribution. The choice of tempering procedure has to be based on the model (or, at least, on the toy-model) leading to this distribution to preserve practical sense of the model, and, in the same time, to make more realistic this ideal model. Unfortunately, there are almost no toy-models in economics to explain the appearance of heavy tailed distributions. It would be nice to have such models to understand what type of tempering should be applied.

\section*{Acknowledgements}
The author are grateful to professors Svetlozar T. Rachev and Frank J. Fabozzi for their helpful comments and advices.


\begin{thebibliography}{500}

\bibitem{M1963}
\newblock Benoit Mandelbrot (1963),
\newblock The Variation of Certain Speculative Prices
\newblock The Journal of Business, Vol. 36, No. 4, pp. 394-419

\bibitem{F1965}
\newblock Fama, E. F. (1965)
\newblock The behavior of stock market prices
\newblock Journal of Business, vol. 38, 34-105.

\bibitem{S1967}
\newblock Samuelson P.A. (1967)
\newblock Efficient Portfolio Selection for Pareto-Levy Investments.
\newblock Journal of Financial and Quantitative Analysis, 2:107-122.

\bibitem{Em1997}
\newblock EMBRECHTS, P., KLUPPELBERG, C., AND MIKOSCH, T. (1997)
\newblock Modelling  Extremal  Events  for  Insurance  and  Finance.
\newblock Berlin:Springer-Verlag

\bibitem{RM2000}
\newblock Svetlozar T. Rachev, Stefan Mittnik (2000)
\newblock Stable Paretian Models in Finance,
\newblock Wiley \& Sons, New York.


\bibitem{KRSz}
\newblock Klebanov L.B., Rachev S.T., Szekely G.J. (1999)
\newblock Pre-limit Theorems and Their Applications, 
\newblock Acta Applicandae Mathematicae,  v.58, 159-174.

\bibitem{KRSa}
\newblock  Klebanov L.B., Rachev S.T., Safarian M. (2000)
\newblock Local Pre-Limit Theorems and Their Applications to Finance, 
\newblock Applied Mathematics Letters, v. 3, 2000, 73-78.

\bibitem{CGMY}
\newblock P. Carr, H. Geman, D. Madan, M. Yor (2002) 
\newblock The Fine Structure of Asset Returns: An Empirical Investigation, 
\newblock Journal of Business, 75 (2), 305–332.

\bibitem{Fel2}
\newblock William Feller (1971)
\newblock An Introduction to Probability Theory and Its Applications, vol. 2,
\newblock John Wiley \& Sons, New York

\bibitem{ST}
\newblock G. Samorodnitsky and M.S. Taqqu (1994)
\newblock Stable Non-Gaussian Processes: Stochastic Models with Infinite Variance,
\newblock Chapman \& Hall.

\bibitem{KMR}
\newblock Klebanov L.B.,  Melamed J.A.,  Rachev S.T. (1989)
\newblock On the products of a random number of random variables in connection
with a problem from mathematical economics, 
\newblock Lecture Notes in Mathematics (Springer-Verlag),  Vol.1412, 103-109.

\bibitem{BL}
\newblock S. I. Boyarchenko, S. Z. Levendorskiǐ (2000)
\newblock  Option pricing for truncated Lévy processes,
\newblock International Journal of Theoretical and Applied Finance, 3 (3), 549–552.

\bibitem{CGMY}
\newblock P. Carr, H. Geman, D. Madan, M. Yor (2002) 
\newblock The Fine Structure of Asset Returns: An Empirical Investigation, 
\newblock Journal of Business, 75 (2), 305–332.

\bibitem{KRBF}
\newblock Kim, Y.S.; Rachev, Svetlozar. T.;, Bianchi, M.L.; Fabozzi, F.J. (2007)
\newblock A New Tempered Stable Distribution and Its Application to Finance. 
\newblock In: Georg Bol, Svetlozar T. Rachev, and Reinold Wuerth (Eds.), Risk Assessment: Decisions in Banking and Finance, Physika Verlag, Springer.

\bibitem{SK}
\newblock Lenka Sl\'{a}mov\'{a} and Lev B. Klebanov (2014)
\newblock On discrete approximations of stable distributions,
\newblock arXiv:1403.3671 [math.PR], 1-7.


\end{thebibliography}
\end{document}